\documentclass[10pt]{amsart}

\usepackage[utf8]{inputenc}
\usepackage[T1]{fontenc}
\usepackage[english]{babel}
\usepackage{graphicx}
\usepackage{amsmath}
\usepackage{amssymb}
\usepackage{amsthm}
\usepackage{fullpage}
\usepackage[left=2.4cm, right=2.4cm, top=1.2cm, bottom=2.3cm]{geometry}
\usepackage{mathrsfs}
\usepackage{mathtools}
\usepackage[hidelinks]{hyperref}
\usepackage{amsaddr}
\usepackage{cleveref}

\hypersetup{allcolors=black}

\newtheorem{theorem}{Theorem}
\newtheorem{definition}{Definition}

\newcommand\N{\mathbb N}

\newcommand\R{\mathbb R}

\newcommand\Cc{\mathcal C}

\author{Pegah Pournajafi}
\address{Chaire Combinatoire, Collège de France, Université PSL, 75005, Paris, France}

\title{A note on $\chi$-unbounded classes of geometric graphs}

\begin{document}

\begin{abstract}
	{\small
		We show that there exist infinitely many classes of intersection graphs of geometric objects that are not $\chi$-bounded -- namely, $d$-CBU graphs for $d\geq 3$ -- and each is incomparable with the class of Burling graphs. This answers a folklore open problem on whether Burling graphs are the sole source of unbounded chromatic number among geometric intersection classes.
	}
\end{abstract}

\maketitle

This short note addresses the following question: Are there classes of intersection graphs of geometric objects that are not $\chi$-bounded and do not contain Burling graphs?
In 2023, a positive answer to this question was given in my PhD thesis~\cite[Appendix B]{PournajafiPhDThesis}.
Recent discussions have brought to my attention that this question is of broader interest, motivating this short note to make the result more accessible to the community.

\smallskip

Let $\mathcal{S} = \{S_1, \dots, S_n\}$ be a collection of sets. The \emph{intersection graph of $\mathcal{S}$} is the graph whose vertex set is~$\mathcal S$ and whose edge set is $\{S_iS_j \mid i\neq j, S_i \cap S_j \neq \varnothing \}$. When each~$S_i$ is a subset of $\mathbb{R}^n$, we call~$G$ a \emph{geometric intersection graph}, or simply a \emph{geometric graph}.
A hereditary class $\Cc$ of graphs is called \emph{$\chi$-bounded} if there exists a function $f \colon \N \to \N $ such that for every~$G \in \Cc$ we have $\chi(G) \leq f(\omega(G))$, where $\chi(G)$ and $\omega(G)$ denote the chromatic number and the clique number of $G$ respectively. If a class is not $\chi$-bounded, we say that it is \emph{$\chi$-unbounded}.

A well-known example of a $\chi$-unbounded geometric graph class is the class of \emph{Burling graphs}, first introduced by Burling in 1965~\cite{Burling1965} as intersection graphs of axis-parallel boxes in $\R^3$. It was later shown to be (strictly) contained in the class of intersection graphs of line-segments \cite{Pawlik2014linesegment} and inhomogeneous homothety-translations of any compact path-connected subset of $\R^2$ \cite{pawlik2013general}.  In~\cite{Pournajafi2022BGasIG}, the author gives the exact characterisation of Burling graphs as a precise subclass of intersection graphs of compact path-connected subsets of $\R^2$, as well as in the case where the sets are all inhomogeneous homothety-translations of the same set. The special cases where the sets are either line-segments or \emph{frames} (i.e., the borders of axis-aligned rectangles with non-empty interior) were introduced earlier in~\cite{BG1PournajafiTrotignon2021} by Trotignon and the author. In this note, we use the definition with frames which we briefly recall here (for a more detailed treatment, see \cite[Definition 6.2.]{BG1PournajafiTrotignon2021} and \cite[Theorem~10]{Pournajafi2022BGasIG}).

\begin{definition}
	A Burling graph is a triangle-free intersection graph of a set of frames in $\R^2$ satisfying: 
	\begin{itemize}
		\item the left side of any frame does not intersect any other frame, and if the right side of a frame intersects another frame, it intersects both its top and bottom sides,
		\item if two frames intersect, then no third frame is contained in the intersection of the regions bounded by the two frames, 
		\item if one frame lies entirely inside the interior of another, and if they both intersect a third frame $F$, then only the top and bottom sides (and not the right side) of $F$ are intersected by the two other frames.
	\end{itemize}
\end{definition}

A key feature of Burling graphs is that many known $\chi$-unbounded geometric graph classes contain Burling graphs, and this often plays a crucial for proving their $\chi$-unboundedness. Examples include line segment graphs, restricted frame graphs, and others (see~\cite{Chalopin2016, Pawlik2014linesegment, pawlik2013general, Krawczyk2014}). Moreover, it has been recently claimed that any strict subclass of Burling graphs is $\chi$-bounded (see~\cite{RzazewskiWalczak2024} for a reference to this work under preparation).

Such observations have prompted a natural question: Are Burling graphs the only geometric graph class causing $\chi$-unboundedness? In other words, do all $\chi$-unbounded geometric graph classes contain Burling graphs? The aim of this note is to answer this question. 
\begin{theorem} \label{thm}
	There exist infinitely many triangle-free $\chi$-unbounded classes of geometric intersection graphs (namely, $d$-CBU graphs for $d\geq 3$) that neither contain nor are contained in the class of Burling graphs. 
\end{theorem}

The class of \emph{Contact graphs of Boxes with Unidirectional contacts}, abbreviated as \emph{CBU graphs}, was first introduced in dimension three by Magnant and Martin in 2011~\cite{MagnantMartin2011}. This was later generalized to higher dimensions by Gonçalves, Limouzy, and Ochem in 2023~\cite{Goncalves2023}. Let us recall their definition.

\begin{definition}
	Let $d\geq 1$ be an integer. A graph is a \emph{$d$-CBU graph} if it is the intersection graph of a family of axis-parallel boxes in $\R^d$ such that for every pair of distinct boxes, their intersection lies in a hyperplane perpendicular to $e_1$, where $e_1, \dots, e_d$ is the standard basis of $\R^d$.
	
	A graph is a \emph{CBU graph} if it is a $d$-CBU graph for some $d \geq 1$.
\end{definition}

Gonçalves, Limouzy, and Ochem show that for every $d\geq 1$, the class of $d$-CBU graphs is triangle-free \cite[Claim 1]{Goncalves2023} and a strict subclass of {$(d+1)$-CBU} graphs \cite[Theorem 19]{Goncalves2023}. Note that $2$-CBU graphs (and therefore $1$-CBU graphs) are subclasses of intersection graphs of axis-parallel boxes in $\R^2$, and thus are $\chi$-bounded by a result of Asplund and Gr\"unbaum \cite[Section~3]{Asplund1960}. Magnant and Martin \cite[Theorem~3]{MagnantMartin2011} showed that the class of 3-CBU graphs is $\chi$-unbounded, thus so is the class of $d$-CBU graphs for $d\geq 3$.

\begin{figure}
	\begin{center}
		\vspace*{-1.6cm}
		\includegraphics[width=15.3cm]{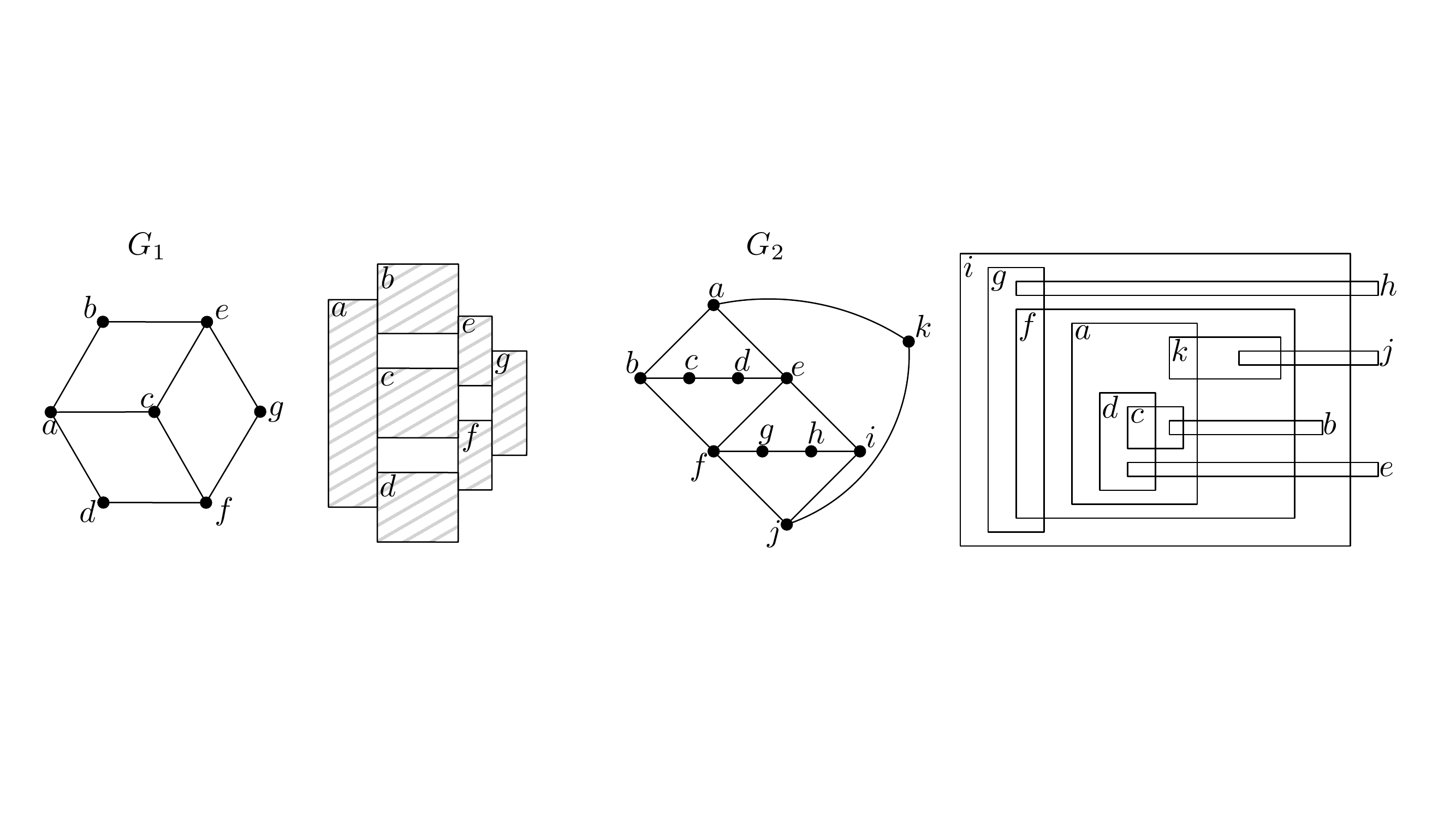}
		\vspace*{-3cm}
		\caption{The graphs $G_1$ and $G_2$.} \label{bothfigs}
	\end{center}
\end{figure} 

\smallskip
We now prove~\Cref{thm}. For a slightly different proof see~\cite[Appendix B]{PournajafiPhDThesis}.
\begin{proof}[Proof of~\Cref{thm}]
	Let $d\geq 3$. We show that the class of $d$-CBU graphs satisfies the desired properties. Recall that $d$-CBU graphs are triangle-free and $\chi$-unbounded.
	
	On the one hand, as shown in~\Cref{bothfigs}, the graph $G_1$ is a 2-CBU graph, hence also a $d$-CBU graph. It was proved by Trotignon and the author in~\cite[Theorem 24]{Pournajafi2020M2} that wheels are not Burling graphs (see also~\cite{BG2PournajafiTrotignon2021} and~\cite{Davies2021}). Since~$G_1$ is a wheel, it is not a Burling graph. 
	
	On the other hand,~\Cref{bothfigs} represents the graph $G_2$ as a Burling graph. Gon\c{c}alves, Limouzy, and Ochem \cite[Lemma 16]{Goncalves2023} have shown that $G_2$ is not a CBU graph, and thus not a $d$-CBU graph.
	
	Hence, the class of Burling graphs and the class of $d$-CBU graphs are incomparable under inclusion. 
	This shows that the class of $d$-CBU graphs for $d\geq 3$, as well as the class of CBU-graphs, satisfy the desired conditions. The fact that $d$-CBU graphs are two-by-two distinct \cite[Theorem 19]{Goncalves2023} completes the proof.
\end{proof}

\end{document}